
\magnification=1200

\overfullrule=0pt

\input amstex

\documentstyle{amsppt}

\hsize=165truemm

\vsize=227truemm


\def\p#1#2{{{\Bbb P}^{#1}_{#2}}}

\def\im{\operatorname{im}}

\def\rk{\operatorname{rk}}

\def\ext{\operatorname{Ext}}

\def\coker{\operatorname{coker}}

\def\hom{\operatorname{Hom}}

\def\Hom#1{{{\Cal H}\kern -0.25ex{\italic om\/}_{\Ofa#1}}}

\def\Syz{{\Cal S}\kern -0.25ex{\italic yz\/}}

\def\Ofa#1{{{\Cal O}_{#1}}}

\def\Tor{\mathop{{\Cal T}\kern -0.30ex{\italic
or\/}}\nolimits}

\def\Ext#1#2{{{\Cal E}\kern -0.25ex{\italic
xt\/}^#1_{\Ofa{#2}}}}

\def\mapright#1{\mathbin{\smash{\mathop{\longrightarrow}
\limits^{#1}}}}

\def\mapdown#1{\Big\downarrow\rlap{$\vcenter{\hbox
{$\scriptstyle#1$}}$}}


\topmatter

\title
Quadratic sheaves and self--linkage
\endtitle

\author
G. Casnati and F. Catanese
\endauthor

\address
\language=2
Gianfranco Casnati: Dipartimento di Matematica, Politecnico di Torino, c.so
Duca degli Abruzzi 24, I--10129 Torino (Italy)
\language=0
\endaddress

\email
casnati\@calvino.polito.it
\endemail

\address
\language=2
Fabrizio Catanese: Mathematisches Institut der Universit\"at, Bunsenstrasse
3--5,D--37073 G\"ottingen (Germany)
\language=0
\endaddress

\email
catanese\@cfgauss.uni-math.gwdg.de
\endemail



\endtopmatter

\document

\head
0. Introduction and first results
\endhead

The present paper is devoted to the proof of a structure theorem for
self--linked pure subschemes
$C\subseteq\p nk$ of codimension $2$ over a field $k$ of characteristic $p\ne2$.

We use the symbol $C$ because the classical case
to be studied was the one of reduced curves in $\p3k$. In this case, one
can describe
the notion of self--linkage in non--technical terms, saying that $C$ is
self--linked if
and only if there are surfaces $F$ and $G$ such that their complete
intersection is the
curve $C$ counted with multiplicity $2$. This is a special case of the
notion of
linkage ($C$ is linked to $C'$ if $C\cup C'$ is the complete intersection
of two
surfaces $F$ and $G$), classically introduced by R. Apery, F. Gaeta (see
[Ap], [Gae]) and
later deeply investigated by the algebraic point of view by Ch. Peskine and
L. Szpiro
and by P. Rao (see [P--S], [Rao1]). The special case of self--linkage was
however studied
before, in the work of E. Togliatti (see [To1], [To2] ), and later D.
Gallarati (see
[Gal]), in the form of the theory of contact between surfaces.

In [Ca1] the theory of
contact was related to a new theory, of the so called even sets of nodes,
and later Rao
used these ideas to obtain a structure theorem for projectively Cohen--Macaulay
self--linked subschemes of codimension $2$ in projective spaces (see [Rao2]).

Recently, Walter's structure theorem (see [Wa]) for subcanonical subschemes of
codimension
$3$ opened the way to solving some old conjectures about even sets of nodes
and contact
of hypersurfaces (see [C--C]).

A basic ingredient was the algebraic concept of quadratic sheaves,
generalizing to a
greater extent the geometric notion of contact and even sets. This notion
was applied in
[C--C] to the classification of even sets of nodes on a surface
$F\subseteq\p3k$
for low values of the degree $d$.

On the other hand, let $F\subseteq\p3k$ be a surface whose only
singularities are an even
set of nodes $\Delta$. Then there is a curve $C$ on $F$ passing through
the points of $\Delta$, and a surface $G$ such that $F\cap G=2C$ as cycles
(see [Ca1], [Gal]).
With this in mind it is therefore only natural to apply the structure
theorem for
quadratic sheaves in order to obtain a structure theorem for self--linkage.

This is done in the present paper, where we generalize the previously cited
result
of P. Rao to arbitrary pure subschemes of $\p nk$ of codimension $2$,
(see also the survey [Ca2] of the second author for a preliminary version
of these
results).

\proclaim{Main Theorem}
Let $k$ be a field of characteristic $p\ne2$ and $C\subseteq\p nk$ be a pure
subscheme of codimension $2$ which is self--linked through hypersurfaces
$F:=\{f=0\}$ and
$G:=\{g=0\}$ of respective degrees $d$, $m$. Let $\Im_C$ be its sheaf of
ideals and set
${\Cal F}_F:=\Im_C/f\Ofa{\p nk}(-d)$, ${\Cal
F}_G:=\Im_C/g\Ofa{\p nk}(-m)$.

Assume moreover that if $n\equiv1\mod 4$ and $d+m-n-1=2\varrho$, then the
following two
equivalent congruences hold:
$$
\aligned
\chi(\Ofa{\p nk}(\varrho))-\chi(\Ofa{\p nk}(\varrho-d))-\chi(\Ofa
C(\varrho))&\equiv0\mod2,\\
\chi(\Ofa{\p nk}(\varrho))-\chi(\Ofa{\p nk}(\varrho-m))-\chi(\Ofa
C(\varrho))&\equiv0\mod2.
\endaligned
$$

Then there exist a locally free $\Ofa{\p nk}$--sheaf
${\Cal E}$, a symmetric map $\alpha\colon \check{\Cal
E}(-d-m)\to{\Cal E}$ and a resolution
$$
0 \longrightarrow\Ofa{\p nk}(-d)\oplus\check{\Cal
E}(-d-m)\mapright{{{\gamma\ {}^t\lambda}\choose {\lambda\ \ \alpha}}}
\Ofa{\p nk}(-m)\oplus{\Cal E}
\longrightarrow{\Cal F}_F\longrightarrow0\tag 1
$$
inducing exact sequences
$$
\gather
0 \longrightarrow\check{\Cal E}(-d-m)\mapright{{{}^t\lambda}\choose
{\alpha}}
\Ofa{\p nk}(-m)\oplus{\Cal E}
\longrightarrow\Im_C\longrightarrow0,\tag 2\\
0 \longrightarrow\check{\Cal E}(-d-m)\mapright{
{\alpha}}
{\Cal E}
\longrightarrow{\Cal F}_G\longrightarrow0.\tag 3
\endgather
$$

Conversely, given a subscheme $C$ of codimension $2$, assume that there does
exist a sequence (2) with the above property of
$\alpha$ being symmetric and $\deg(\det(\alpha))=m$. Then
$C$ is self--linked through the hypersurfaces
$F$ and
$G$ of respective equations
$f:=\det\pmatrix0&{}^t\lambda\\
\lambda&\alpha\endpmatrix$ and $g:=\det(\alpha)$.
\qed
\endproclaim

Moreover, in section 3 we discuss the problem whether the self--linkage ideal
$(f,g)$ is uniquely determined once the generator $g$ of higher degree is
fixed, and prove that the
answer is positive  under the condition that the subscheme be locally
Gorenstein.

\subhead
Acknowledgements
\endsubhead
Both the authors acknowledge support from the projects AGE, EAGER and Vigoni.

\subhead
First results
\endsubhead
For the reader's benefit we recall the following definitions and results
proved from
[C--C]. From now on we always assume that $k$ is a field of characteristic
$p\ne2$.

\definition{Definition 0.1}
Let $X$ be a projective, locally Cohen--Macaulay scheme. We say that a
coherent, locally Cohen--Macaulay  sheaf of $\Ofa X$--modules $\Cal F$ is a
$\delta/2$--quadratic sheaf on $X$, $\delta=0,1$, if
there exists a symmetric bilinear map
$$
{\Cal F}\times{\Cal F}\to\Ofa X(-\delta)
$$
inducing an isomorphism $\sigma\colon{\Cal
F}(\delta)\mapright\sim\Hom{X}\big({\Cal F},\Ofa X\big)$.
\enddefinition

\remark{Remark 0.2}
If ${\Cal F}$ is $\delta/2$--quadratic, then it is reflexive since the natural
map ${\Cal F}\to{\Cal F}\check{\ }\check{\ }$ equals
$\check\sigma^{-1}\circ\sigma(-\delta)$.
\endremark
\medbreak

Assume now that $F\subseteq\p nk$ is a hypersurface of degree $d$
and that
$\Cal F$ is a $\delta/2$--quadratic sheaf on $F$.

The main result of sections 1 and 2 of [C--C] concerns
a characterization of quadratic sheaves on
hypersurfaces in $\p nk$ (including the needed parity condition as pointed out
in theorem 9.1
of [E--P--W]).

\proclaim{Theorem 0.3}
Let $F\subseteq\p nk$ be a hypersurface of degree $d$ and let
$\Cal F$ be a
$\delta/2$--quadratic sheaf on $F$. Then
$\Cal F$ fits into an exact sequence of the form
$$
0\longrightarrow\check{\Cal E}(-d-\delta)\mapright\varphi{\Cal
E}\longrightarrow{\Cal F}\longrightarrow0\tag0.3.1
$$
where $\Cal E$ is a locally free $\Ofa{\p nk}$--sheaf and $\varphi$ is a
symmetric map
if and only if the following parity condition holds: if $n\equiv1\mod 4$ and
$n+1-d-\delta=2r$, then also $\chi({\Cal F}(-r))$ is even.

Moreover we can choose $\Cal E$ such that $H^j_*\big(\p nk,{\Cal E}\big)=0$ for
$n>j>(n-1)/2$.
\qed
\endproclaim

\remark{Remark 0.4}
Following [Ca1] we say that $\Cal F$ is split symmetric if one can choose
$\Cal E$ to be a
direct sum of line bundles. This is possible if and only if $\Cal F$  is
arithmetically Cohen-Macaulay i.e., if and only if $H^i_*\big(\p nk,{\Cal
F}\big) = 0$
for each $i=1,\dots,n-2$.
\endremark
\medbreak

\head
1. From self--linkage to quadratic sheaves
\endhead

Our first main theorem is an application of the theory of quadratic
sheaves. For another
application see [C--C].

\definition{Definition 1.1}
Let $C\subseteq\p nk$ be a pure subscheme of codimension $2$ and let
$\Im_C\subseteq\Ofa{\p nk}$ be its sheaf of ideals:
$C$ is said to be self-linked with respect to the complete intersection
$X:=F\cap G$ of the two hypersurfaces $F$, $G$ of respective degrees $d$,
$m$, if
$C\subseteq X$ and one of the following equivalent conditions holds
\roster
\item"i)" $\Im_X\colon\Im_C=\Im_C$;
\item"ii)" $\Im_C/\Im_X=\Hom{{\p nk}}\big(\Ofa{C},\Ofa{X}\big)$.
\endroster
\enddefinition

For the above well-known equivalence see e.g. theorem 21.23 of [Ei] (see also
[P--S]).

\remark{Remark 1.2}
Indeed, since
$\Im_C$ and $\Im_{X}$ coincide with
$\Ofa{\p nk}$  in codimension $1$ and are torsion free, every
$\phi\in\Hom{{\p nk}}
\big( \Im_C,\Im_{X}\big)$ is given by a rational function, which is in turn
regular by
the normality of $\Ofa{\p nk}$.

Condition i) of definition 1.1 can be thus rewritten as (see [Ca2],
proposition 2.6)
$$
\Im_C= \Hom{{\p nk}}\big(\Im_C,\Im_{X}\big).\tag 1.2.1
$$

By duality for finite maps (see [Ha], exercise II.6.10),
$$\Hom{{\p nk}}\big(\Ofa C,\Ofa{X}\big)=\Hom{{\p nk}}\big(
\Ofa C,\omega_{X}(n+1-m-d)\big)=\omega_C(n+1-m-d).\tag 1.2.2
$$

Finally, let $C$, $F:=\{f=0\}$, $G:=\{g=0\}$ be as above, and assume
$m=\deg(G) \geq
d =\deg(F)$. Then we can replace
$g$ by $g +a f$ ($a$ is a homogeneous polynomial of degree $m-d$), and
obtain by Bertini's theorem that $G$ is smooth outside $C$. Indeed, if $C$
is reduced we
can even have that $G$ is smooth at the generic points of $C$, so that $G$
is a normal,
whence irreducible, hypersurface. Moreover in this case (cfr. [Ca1],
proposition
2.6) condition i) is equivalent to
$cycle(X)=2C$.
\endremark
\medbreak

\definition{Definition 1.3}
Let ${\Cal F}_G:=\Im_C/ g\Ofa{\p nk}(-m)$, ${\Cal F}_F:=\Im_C/ f\Ofa{\p
nk}(-d)$.
\enddefinition

By the above definition it follows that both ${\Cal F}_G$ and ${\Cal F}_F$
are obviously locally Cohen--Macaulay. By
remark 1.2, since
$$
\Im_{X}/g\Ofa{\p nk}(-m)\cong (f\Ofa{\p nk}(-d)+g\Ofa{\p
nk}(-m))/g\Ofa{\p nk}(-m)\cong f\Ofa{G}(-d),
$$
we have pairings
$$
{\Cal F}_G\times{\Cal F}_G\longrightarrow\Im_{X}/g\Ofa{\p nk}(-m)\cong
f\Ofa{G}(-d),
$$
and analogously
$$
{\Cal F}_F \times{\Cal F}_F  \longrightarrow g\Ofa{F} (-m).
$$
To verify that a twist of ${\Cal F}_G$  yields a quadratic sheaf on $G$
(similarly for ${\Cal F}_F$ on $F$) we have to show that the pairing is
perfect, i.e.  ${\Cal F}_G\cong\Hom{{\p nk}}\big({\Cal F}_G, \Ofa{G}
(-d)\big)$.

Let $\psi\in\Hom{{\p nk}}\big({\Cal F}_G,\Ofa{G}(-d)\big)$, then $\psi$ induces
$\psi':\Im_C\to\Ofa G(-d)\cong\Im_{X}/ g\Ofa{\p nk}(-m)$. If $\psi'$ is in
the image of
a $\phi\in\Hom{{\p nk}}\big(\Im_C,\Im_{X}\big)$, then, by remark 1.2,
$\phi\in\Im_C$, and clearly $\phi$  induces a zero $\psi$ if and only if
$\phi\in
g\Ofa{\p nk}(-m)$.

Thus, we are left with the verification that
$$
\Hom{{\p nk}}\big(\Im_C,\Im_{X}\big) \to
\Hom{{\p nk}}\big(\Im_C,\Im_{X}/g\Ofa{\p nk}(-m)\big)
$$
is surjective. By the $\Ext{i}{\p nk}$--exact sequence, the above
surjectivity is
equivalent to the injectivity of
$j\colon\Ext{1}{\p nk}\big(\Im_C, g\Ofa{\p nk}(-m)\big)
\to\Ext{1}{\p nk}\big(\Im_C,\Im_{X}\big)$.

To this purpose, consider
$$
\align
0\longrightarrow\Ext{1}{\p nk}\big(\Ofa{C},
\Im_X\big)&\mapright{i}\Ext{1}{\p nk}\big(\Ofa{C},
\Ofa{G}(-d)\big)\longrightarrow\\
&\longrightarrow\Ext{2}{\p nk}\big(\Ofa{C},
g\Ofa{\p nk}(-m)\big)
\mapright{j'}\Ext{2}{\p nk}\big(\Ofa{C},
\Im_X\big)
\endalign
$$
where $j'$ corresponds to $j$ via the natural equivalence
$\Ext{1}{\p nk}\big(\Im_{C},
\cdot\big)\cong\Ext{2}{\p nk}\big(\Ofa{C},
\cdot\big)$.
It follows that $j$ is injective if and only if $i$ is an isomorphism.
We have the exact sequence
$$
\align
\Hom{{\p nk}}\big(\Ofa{\p nk},
\Im_X\big)&\longrightarrow\Hom{{\p nk}}\big(\Im_{C},
\Im_X\big)\longrightarrow\\
&\longrightarrow\Ext{1}{\p nk}\big(\Ofa{C},
\Im_X\big)\longrightarrow\Ext{1}{\p nk}\big(\Ofa{\p nk},
\Im_X\big)\cong0.
\endalign
$$
Since $\Hom{{\p nk}}\big(\Ofa{\p nk},
\Im_X\big)\cong\Im_X$ and $\Hom{{\p nk}}\big(\Im_{C},
\Im_X\big)\cong\Im_C$ then
$$
\Ext{1}{\p nk}\big(\Ofa{C},\Im_X\big)\cong\Im_C/\Im_X\cong\Hom{{\p
nk}}\big(\Ofa{C},
\Ofa{X}\big)\cong\omega_C(n+1-m-d).
$$
We also have
$$
\align
\Ext{1}{\p nk}\big(\Ofa{C},
\Ofa{\p nk}(-d)\big)&\longrightarrow\Ext{1}{\p nk}\big(\Ofa{C},
\Ofa{G}(-d)\big)\longrightarrow\\
&\longrightarrow\Ext{2}{\p nk}\big(\Ofa{C},
g\Ofa{\p nk}(-m-d)\big)\mapright{g}\Ext{2}{\p nk}\big(\Ofa{C},
\Ofa{\p nk}(-d)\big).
\endalign
$$
Since $C\subseteq G$, then the multiplication by $g$ is zero. Moreover,
$\Ext{1}{\p nk}\big(\Ofa{C},
\Ofa{\p nk}(-d)\big)=0$ since $C$ is locally Cohen--Macaulay. We obtain that
$$
\Ext{1}{\p nk}\big(\Ofa{C},
\Ofa{G}(-d)\big)\cong\Ext{2}{\p nk}\big(\Ofa{C},
g\Ofa{\p nk}(-m-d)\big)\cong\omega_C(n+1-m-d).
$$
The following easy lemma thus concludes the proof that
$i\colon\omega_C(m+d-n-1)
\to\omega_C(m+d-n-1)$ is an isomorphism.

\proclaim{Lemma 1.4}
If $\Cal H$ is a coherent $\Ofa{\p nk}$--sheaf and $\varrho\colon{\Cal
H}\to{\Cal
H}$ is either injective or surjective, then it is an isomorphism.
\endproclaim
\demo{Proof}
Let $\varrho$ be injective (resp. surjective) and ${\Cal
K}:=\coker(\varrho)$ (resp.
${\Cal K}:=\ker(\varrho)$). By Serre'stheorem B we get
$h^1\big(\p nk,{\Cal H}(t)\big)=0$ (resp. $h^1\big(\p nk,{\Cal K}(t)\big)=0$)
for each $t$ large enough, hence
$h^0\big(\p nk,{\Cal K}(t)\big)=0$ for $t$ large enough, whence ${\Cal K}=0$
by Serre's theorem A.
\qed
\enddemo

Set $d = 2d' + \delta$, $m = 2 m' + \mu$, $\delta,\mu\in\{0,1\}$. Then
${\Cal F}_G (d')$ is a $\delta/2$--quadratic sheaf on $G$, and similarly ${\Cal
F}_F (m')$ is a
$\mu/2$--quadratic sheaf on $F$. According to theorem 0.3, we obtain two
locally
free
$\Ofa{\p nk}$--sheaves ${\Cal E}'_G$ and ${\Cal E}'_F$ provided that the
respective
parity conditions hold.

It is convenient to rewrite such conditions. If $n\equiv1\mod 4$ and
$d+m-n-1=2\varrho$, then
we want that $\chi({\Cal
F}_F(\varrho))\equiv\chi({\Cal F}_G(\varrho))\equiv0\mod2$. By
two obvious exact sequences the above conditions are equivalent to the two
congruences
$$
\aligned
\chi(\Ofa{\p nk}(\varrho))-\chi(\Ofa{\p nk}(\varrho-d))-\chi(\Ofa
C(\varrho))&\equiv0\mod2,\\
\chi(\Ofa{\p nk}(\varrho))-\chi(\Ofa{\p nk}(\varrho-m))-\chi(\Ofa
C(\varrho))&\equiv0\mod2.
\endaligned\tag 1.5
$$
In turn the two congruences above are equivalent each other. In fact, it
suffices to show that
$$
\chi(\Ofa{\p nk}(\varrho-d))+\chi(\Ofa{\p
nk}(\varrho-m))=0.
$$

To this purpose recall that
$$
\chi(\Ofa{\p nk}(h))={{(h+n)\dots(h+1)}\over{n!}}.
$$
In our case $n=4a+1$ and $\varrho$ is even. Since, by the linkage
condition, $dm$ is even too, the parity of $d+m-n-1$
yields
$d=2d'$,
$m=2m'$, hence
$$
\chi(\Ofa{\p
nk}(\varrho-d))={{(2a+i)\dots(-2a+i)}\over{(4a+1)!}}=-{{(2a-i)\dots(-2a-i)}\over
{(4a+1)!}}=\chi(\Ofa{\p
nk}(\varrho-m)),
$$
where $i=m'-d'$.

We define
${\Cal E}_G:={\Cal E}'_G(-d')$ and ${\Cal E}_F:={\Cal E}'_F(-m')$. Then theorem
0.3 rewrites as follows.

\proclaim{Proposition 1.6}
Assume that if $n\equiv1\mod 4$ and
$n+1-d-m=2\varrho$ then the two equivalent congruences (1.5) above hold.

Then there exist two locally free $\Ofa{\p nk}$--sheaves ${\Cal
E}_G$ and
${\Cal E}_F$, such that $H^i_*\big(\p nk,{\Cal E}_G\big)=H^i_*\big(\p nk,{\Cal
E}_F\big)=0$ for $n>i>(n-1)/2$, fitting into exact sequences
$$
\aligned
&0\longrightarrow
\check{\Cal E}_F(-d-m)\mapright{\alpha_F}{\Cal E}_F\longrightarrow{\Cal F}_F
\longrightarrow0,\\
&0\longrightarrow
\check{\Cal E}_G(-d-m)\mapright{\alpha_G}{\Cal E}_G\longrightarrow{\Cal F}_G
\longrightarrow0,
\endaligned\tag1.6.1
$$
where ${}^t\alpha_F=\alpha_F$,
${}^t\alpha_G=\alpha_G$.
\qed
\endproclaim

In the above proposition and in what follows, the superscript $t$ denotes the
dual morphism twisted by $-d-m$.

\remark{Remark 1.7}
Let $\ell\subseteq\p nk$ be a line disjoint from $C$. Then the restrictions of
the sequences (1.6.1) to $\ell$ are still exact. On the other hand ${\Cal
F}_F{}_{\vert
\ell}\cong\Ofa {F\cap \ell}$ and ${\Cal F}_G{}_{\vert
\ell}\cong\Ofa {G\cap \ell}$, hence
$$
c_1({\Cal E}_F)=-{1\over2}\rk({\Cal E}_F)(d+m)+{d\over2},\qquad c_1({\Cal
E}_G)=-{1\over2}\rk({\Cal E}_G)(d+m)+{m\over2}.
$$
\endremark
\medbreak

\head
2. From quadratic sheaves to self--linkage
\endhead

From now on we will assume that the parity condition holds.

Let $\Ofa{\p
nk}(-d)\mapright{f}\Ofa{\p nk}$ be the multiplication by $f$. Then it induces a
map $\vartheta\colon\Ofa{\p nk}(-d)\to{\Cal F}_G$. Since we have a chain of
homomorphisms
$$
\aligned
\hom_{\Ofa{\p nk}}\big(\Ofa{\p nk}(-d),{\Cal
E}_G\big)&\longrightarrow\hom_{\Ofa{\p nk}}\big(\Ofa{\p
nk}(-d),{\Cal
F}_G\big)\longrightarrow\\
&\mapright\partial\ext^1_{\Ofa{\p nk}}\big(\Ofa{\p nk}(-d),\check{\Cal
E}_G(-d-m)\big)\cong\\
&\cong
H^1\big(\p nk,\check{\Cal E}_G(-m)\big)\cong
H^{n-1}\big(\p nk,{\Cal E}_G(m+d-n-1)\big)\check{}=0
\endaligned\tag 2.1
$$
we can lift $\vartheta$ to a map $\lambda\colon \Ofa{\p nk}(-d)\to{\Cal E}_G$.

Argueing as in (2.1) we obtain that the natural map
$$
\hom_{\Ofa{\p nk}}\big({\Cal
E}_G,\Im_C\big)\to\hom_{\Ofa{\p nk}}\big({\Cal E}_G,{\Cal
F}_G\big)
$$
is surjective
whence we infer that the surjection ${\Cal E}_G\to{\Cal F}_G$ can be
lifted to $\nu\colon{\Cal E}_G\to\Im_C$.
Notice that the map
$\nu\circ\lambda$ is given by $\widetilde{f}\in H^0\big(\p nk,\Ofa{\p
nk}(d)\big)$
congruent $f\mod g$.

Let $r:=\rk({\Cal E}_G)$. We have a map
$$
\Lambda^{r-1}(\alpha_G)\colon\Lambda^{r-1}(\check{\Cal
E}_G(-d-m))\to\Lambda^{r-1}{\Cal E}_G.
$$
Since $\Lambda^{r-1}\check{\Cal E}_G\cong{\Cal E}_G\otimes\det({\Cal
E}_G)^{-1}$ and $\Lambda^{r-1}{\Cal E}_G\cong\check{\Cal E}_G\otimes\det({\Cal
E}_G)$, twisting by $\det({\Cal E}_G)^{-1}(-d)$ and taking remark 1.7 into
account,
we obtain a map $\alpha_G^{adj}\colon{\Cal E}_G\to\check{\Cal E}_G(-d)$
such that
$\alpha_G\circ\alpha_G^{adj}(-m)\colon{\Cal E}_G(-m)\to{\Cal
E}_G$ is the multiplication by $g$ and
$\alpha^{adj}_G(d)=\check\alpha^{adj}_G$ since
$\alpha_G$, hence $\alpha^{adj}_G$, is symmetric.

We can
then define
$$
\mu:=\check\lambda(-d)\circ\alpha_G^{adj}\colon{\Cal E}_G\to\Ofa{\p nk}.
$$

\proclaim{Proposition 2.2}
$\Im_C$ coincides with the sheaf of ideals of $\Ofa{\p nk}$ generated by $g$
and $\im(\mu)$.
\endproclaim
\demo{Proof}
Let us consider the second sequence in (1.6.1). We have a chain map
$$
\matrix
0&\hskip-3truemm\longrightarrow&\hskip-3truemm\check{\Cal
E}_G(-d-m)&\hskip-3truemm\mapright{\alpha_G}&\hskip-3truemm{\Cal
E}_G&\hskip-3truemm\longrightarrow&\hskip-3truemm{\Cal F}_G&\hskip-3truemm
\longrightarrow&\hskip-3truemm0&\hskip-3truemm\\
&\hskip-3truemm&\hskip-3truemm\mapdown{\beta}&\hskip-3truemm
&\hskip-3truemm\mapdown{\nu}&\hskip-3truemm&\hskip-3truemm\mapdown{id}
&\hskip-3truemm&\hskip-3truemm&\hskip-3truemm\\
0&\hskip-3truemm\longrightarrow&\hskip-3truemm\Ofa{\p nk}
(-m)&\hskip-3truemm\mapright{g}&\hskip-3truemm\Im_C
&\hskip-3truemm\longrightarrow &\hskip-3truemm{\Cal
F}_G&\hskip-3truemm\longrightarrow&\hskip-3truemm0&\hskip-3truemm
\endmatrix\tag 2.2.1
$$
where $\beta$ is induced by the restriction of $\nu$.

By the mapping cone construction we obtain a resolution
$$
0\longrightarrow\check{\Cal E}_G(-d-m)\mapright{s} \Ofa{\p nk} (-m)\oplus
{\Cal E}_G
\longrightarrow\Im_C\longrightarrow0
$$
where $s$ has components $\beta$, $\alpha_G$. Recall that
$\alpha_G(m)\circ\alpha_G^{adj}$ is the multiplication by $g$ whence
$(\nu\circ\alpha_G)(m)\circ\alpha_G^{adj}=g\nu$, since $(\nu\circ
g)(m)=g\nu$.

Diagram (2.2.1) yields $g\beta=\nu\circ\alpha_G$.
Composing on the right with
$\alpha_G^{adj}\circ\lambda$, we obtain
$$
g\beta(m)\circ\alpha_G^{adj}\circ\lambda=
g\nu\circ\lambda\colon\Ofa{\p nk}(-d)\to\Im_C(m)\subseteq\Ofa{\p nk}(m).
$$
Since $g$ is obviously a non--zero divisor in $\Im_C$, then
$\nu\circ\lambda=\beta(m)\circ\alpha_G^{adj}\circ\lambda\colon\Ofa{\p
nk}(-d)\to\Ofa{\p nk}$ and since $\nu\circ\lambda$ is given by
$\widetilde{f}\equiv f\mod g$ the same is true for
$\beta(m)\circ\alpha_G^{adj}\circ\lambda$.

It follows from the above identities that
$s\circ{}^t\mu=(\nu\circ\lambda(-m),g\lambda(-m))$ is
represented by the product matrix
$$
\pmatrix1&0\\0&\lambda\endpmatrix
\pmatrix \widetilde{f}\\ g\endpmatrix.
$$
Thus we obtain the following commutative diagram
\vglue1pt
$$
\matrix
0&\hskip-3truemm\longrightarrow&\hskip-3truemm\Ofa{\p
nk}(-d-m)&\hskip-3truemm\mapright{{\widetilde{f}\choose
g}}&\hskip-3truemm\Ofa{\p
nk}(-m)\oplus
\Ofa{\p
nk}(-d)&\hskip-3truemm\longrightarrow&\hskip-3truemm\Im_X&\hskip-3truemm
\longrightarrow&\hskip-3truemm0&\hskip-3truemm\\
&\hskip-3truemm&\hskip-3truemm\mapdown{{}^t\mu}&\hskip-3truemm&\hskip-3truemm
\mapdown{
{{1\ 0}\choose
{0\ \lambda}}}&\hskip-3truemm&\hskip-3truemm\Big\downarrow&\hskip-3truemm
&\hskip-3truemm&\hskip-3truemm\\
0&\hskip-3truemm\longrightarrow&\hskip-3truemm\check{\Cal
E}_G(-d-m)&\hskip-3truemm\mapright{s}&\hskip-3truemm\Ofa{\p nk}(-m)\oplus {\Cal
E}_G&\hskip-3truemm\longrightarrow&\hskip-3truemm\Im_C&\hskip-3truemm
\longrightarrow&\hskip-3truemm0&\hskip-3truemm
\endmatrix\tag 2.2.2
$$
whose right column is the inclusion.

The mapping cone of (2.2.2) is a resolution of
$\Im_C/\Im_X\cong\omega_C(n+1-d-m)$
(see (1.2.2)). Therefore the dual of the mapping cone of diagram (2.2.2)
yields a
resolution of
$\Im_C$ (see proposition 2.5 of [P--S]) and $\Im_C$ coincides with the
sheaf of ideals
locally generated by the maximal minors of
$$
\pmatrix\beta&1&0\\
\alpha_G&0&\lambda
\endpmatrix
$$
which is the ideal locally generated by the maximal minors of
$(\alpha_G,\lambda)$. Our
statement follows then from the very definition of
$\mu$ and the identity $\det(\alpha_G)=g$.
\qed
\enddemo

Since $\mu\circ\alpha_G=g({}^t\lambda)$, then $\mu$ induces an endomorphism
$\psi$ of
${\Cal F}_G$ fitting into the following commutative diagram
$$
\matrix
0&\hskip-3truemm\longrightarrow&\hskip-3truemm\check{\Cal
E}_G(-d-m)&\hskip-3truemm\mapright{\alpha_G}&\hskip-3truemm{\Cal
E}_G&\hskip-3truemm\longrightarrow&\hskip-3truemm{\Cal F}_G&\hskip-3truemm
\longrightarrow&\hskip-3truemm0&\hskip-3truemm\\
&\hskip-3truemm&\hskip-3truemm\mapdown{{}^t\lambda}&\hskip-3truemm
&\hskip-3truemm\mapdown{\mu}&\hskip-3truemm&\hskip-3truemm\mapdown{\psi}
&\hskip-3truemm&\hskip-3truemm&\hskip-3truemm\\
0&\hskip-3truemm\longrightarrow&\hskip-3truemm\Ofa{\p nk}
(-m)&\hskip-3truemm\mapright{g}&\hskip-3truemm\Im_C
&\hskip-3truemm\longrightarrow &\hskip-3truemm{\Cal
F}_G&\hskip-3truemm\longrightarrow&\hskip-3truemm0&\hskip-3truemm.
\endmatrix
$$

Since ${\Cal F}_G:=\Im_C/g\Ofa{\p nk}(-m)$, proposition 2.2 implies the
surjectivity of $\psi$. Moreover $\psi$ is also injective by lemma 1.4
above. We can
replace the given surjection $\pi\colon{\Cal E}_G\to{\Cal F}_G$ with
$\psi\circ\pi$,
whence we may also replace $\nu$ by $\mu$ in the arguments of proposition
2.2, thus obtaining the following

\proclaim{Proposition 2.3}
There exists $\gamma\in H^0\big(\p nk,\Ofa{\p nk}(d-m)\big)$ such that
$$
f=\det\pmatrix
\gamma&{}^t\lambda\\
\lambda&\alpha_G
\endpmatrix.
$$
\endproclaim
\demo{Proof}
Recall that we are now assuming $\nu=\mu$, whence by the very
definition of $\mu$, $\mu\circ\lambda\equiv f\mod g$. Since
$$
\mu\circ\lambda=\det\pmatrix
0&{}^t\lambda\\
\lambda&\alpha_G
\endpmatrix,
$$
then we obtain the existence of a $\gamma$ such that
$$
f=\det\pmatrix
0&{}^t\lambda\\
\lambda&\alpha_G
\endpmatrix+\gamma g=\det\pmatrix
0&{}^t\lambda\\
\lambda&\alpha_G
\endpmatrix+\det\pmatrix
\gamma&{}^t\lambda\\
0&\alpha_G
\endpmatrix=\det\pmatrix
\gamma&{}^t\lambda\\
\lambda&\alpha_G
\endpmatrix.\quad\qed
$$
\enddemo

Since $f\ne0$
we get a coherent sheaf $\Cal F$ supported on $F$ and an exact sequence
\vglue1pt
$$
0 \longrightarrow\Ofa{\p nk}(-d)\oplus\check{\Cal
E}_G(-d-m)\mapright{{{\gamma\ {}^t\lambda}\choose {\lambda\ \alpha_G}}}
\Ofa{\p nk}(-m)\oplus{\Cal E}_G
\longrightarrow{\Cal F}\longrightarrow0.
$$
It follows that ${\Cal F}(m')$ is $\mu/2$--quadratic and we can
easily construct the following exact sequence of the vertical complexes $C_i$'s
$$
\matrix
&&C_1&
&C_2&&
C_3&&\\
\\
0&\longrightarrow&\check{\Cal
E}_G(-d-m)&\longrightarrow
&\Ofa{\p nk}(-m)\oplus{\Cal E}_G&
\longrightarrow&\Im_C&\longrightarrow&0\\
&&\mapdown{0\choose {id}}&
&\mapdown{id}&&
\mapdown{\eta}&&\\
0&\longrightarrow&\Ofa{\p nk}(-d)\oplus\check{\Cal
E}_G(-d-m)&\longrightarrow
&\Ofa{\p nk}(-m)\oplus{\Cal E}_G&
\longrightarrow&{\Cal F}&\longrightarrow&0.
\endmatrix
$$

The associated long exact sequence gives $\Ofa{\p nk}(-d)\cong H_1(C_1)\cong
H_0(C_3)\cong\ker(\eta)$. On the other hand $f\in\ker(\eta)$, thus
${\Cal F}\cong{\Cal F}_F$.

If we set ${\Cal E}:={\Cal E}_G$ and
$\alpha:=\alpha_G$, the above discussion proves the \lq\lq only if\rq\rq\
part of
the statement of the following main theorem.

\proclaim{Theorem 2.4}
Let $k$ be a field of characteristic $p\ne2$ and $C\subseteq\p nk$ be a pure
subscheme of codimension $2$ which is self--linked through hypersurfaces
$F:=\{f=0\}$ and
$G:=\{g=0\}$ of respective degrees $d$, $m$. Let $\Im_C$ be its sheaf of
ideals and set
${\Cal F}_F:=\Im_C/f\Ofa{\p nk}(-d)$, ${\Cal
F}_G:=\Im_C/g\Ofa{\p nk}(-m)$.

Assume that if $n\equiv1\mod 4$ and $d+m-n-1=2\varrho$, then the following two
equivalent congruences hold:
$$
\aligned
\chi(\Ofa{\p nk}(\varrho))-\chi(\Ofa{\p nk}(\varrho-d))-\chi(\Ofa
C(\varrho))&\equiv0\mod2,\\
\chi(\Ofa{\p nk}(\varrho))-\chi(\Ofa{\p nk}(\varrho-m))-\chi(\Ofa
C(\varrho))&\equiv0\mod2.
\endaligned
$$

Then there exist a locally free $\Ofa{\p nk}$--sheaf
${\Cal E}$, a symmetric map $\alpha\colon \check{\Cal
E}(-d-m)\to{\Cal E}$ and a resolution
$$
0 \longrightarrow\Ofa{\p nk}(-d)\oplus\check{\Cal
E}(-d-m)\mapright{{{\gamma\ {}^t\lambda}\choose {\lambda\ \ \alpha}}}
\Ofa{\p nk}(-m)\oplus{\Cal E}
\longrightarrow{\Cal F}_F\longrightarrow0\tag 2.4.1
$$
inducing exact sequences
$$
\gather
0 \longrightarrow\check{\Cal E}(-d-m)\mapright{{{}^t\lambda}\choose
{\alpha}}
\Ofa{\p nk}(-m)\oplus{\Cal E}
\longrightarrow\Im_C\longrightarrow0,\tag 2.4.2\\
0 \longrightarrow\check{\Cal E}(-d-m)\mapright{
{\alpha}}
{\Cal E}
\longrightarrow{\Cal F}_G\longrightarrow0.\tag 2.4.3
\endgather
$$

Conversely, given a subscheme $C$ of codimension $2$, assume that there does
exist a sequence (2.4.2) with the above property of
$\alpha$ being symmetric and $\deg(\det(\alpha))=m$. Then
$C$ is self--linked through the hypersurfaces
$F$ and
$G$ of respective equations
$f:=\det\pmatrix0&{}^t\lambda\\
\lambda&\alpha\endpmatrix$ and $g:=\det(\alpha)$.
\endproclaim
\demo{Proof}
There remains to prove the converse assertion. Let $P\in\p nk$ and consider in
$\Ofa{\p nk,P}$ the maximal minors $f_i$ of the matrix $(\lambda\ \alpha)$
obtained by deleting the
$i^{\roman th}$ column: in particular $g=f_1$. Recall that
$\Im_{C,P}=(f_1,\dots,f_{r+1})$ in
$\Ofa{\p nk,P}$, hence for each pair of indices $i,j=2,\dots r+1$ we have
$$
u_{i,j}f+v_{i,j}g=f_if_j
$$
for suitable $u_{i,j},v_{i,j}\in\Ofa{\p nk,P}$ (the determinantal identity (1.2)
of [Ca1] with
$k=j=1$).

It follows that
neither $f$ nor $g$ are identically zero, else $C$ would have
a codimension one component. The same identity shows that
$\Im_C^2\subseteq\Im_X$. Moreover it is always true that
$\Im_X\subseteq\Im_C$, thus $F$ and $G$ have no common components.

We now prove that $\Im_X:\Im_C=\Im_C$. To this purpose we first check that
${\Cal
F}:=\coker(\alpha)\cong\Im_C/g\Ofa{\p nk}(-m)$. Indeed the diagram
$$
\matrix
0&\longrightarrow&\check{\Cal
E}_G(-d-m)&\longrightarrow
&\Ofa{\p nk}(-m)\oplus{\Cal E}_G&
\longrightarrow&\Im_C&\longrightarrow&0\\
&&\mapdown{{id}}&
&\mapdown{(0\ id)}&&
&&\\
0&\longrightarrow&\check{\Cal
E}_G(-d-m)&\longrightarrow
&{\Cal E}_G&
\longrightarrow&{\Cal F}&\longrightarrow&0.
\endmatrix
$$
induces a surjection $\Im_C\to{\Cal F}$. Argueing as above we obtain that
its kernel is $g\Ofa{\p nk}(-m)$.

The map $\alpha$ induces an isomorphism ${\Cal F}\cong\Hom{{\p
nk}}\big({\Cal F}, \Ofa{G} (-d)\big)$, since $\deg(f)=d$. We have
$h\in\Im_{X,P}:\Im_{C,P}$ if
and only if
$h\Im_{C,P}\subseteq\Im_{X,P}$.

Taking residue classes $\mod g$, this is equivalent to $h {\Cal
F}_P\subseteq\Im_{X,P}/g\Ofa{\p nk,P}(-m)\cong\Ofa
{G,P}(-d)$, since $g\in\Im_X$, hence to the fact that $h\mod g$ is in $\Hom{{\p
nk,P}}\big({\Cal F}_P, \Ofa{G,P} (-d)\big)\cong{\Cal F}_P$, i.e.
$h\in\Im_{C,P}$.
\qed
\enddemo

\head
3. Self--linkage of generically Gorenstein subschemes
\endhead

In this section we shall inspect more deeply the case when $C$ is a
generically Gorenstein pure subscheme of codimension $2$ of $\p nk$ which is
self--linked through two hypersurfaces $F:=\{f=0\}$ and $G:=\{g=0\}$ of degrees
$\deg(F)=:d\le m:=\deg(G)$.

The first hypothesys implies that $C$ is also generically locally complete
intersection
since its codimension is $2$. Owing to the isomorphism (1.2.2) we have
$$
\Im_C/\Im_X\cong\omega_C(n+1-m-d)
$$
whence at each generic point $P\in C$ the sheaf $\omega_{C,P}$ is
invertible and a lift of a
generator yields $y$ such that
$\Im_{C,P}=(f,g,y)\Ofa{\p nk,P}$. Since $C$ is generically complete
intersection and
$\Im_X=(f,g)$ then $\Im_C$ is either $(f,y)$ or $(g,y)$ at $P$. In the
first case we have locally $g=af+by$. By changing globally $g$
with $g+cf$ for a suitable $c$ we can assume that $a$ is invertible at each
generic point of $C$, hence at each generic point $P\in C$ we may assume
that we are in the case $\Im_C=(g,y)$ holds generically.

Let now $H:=\{h=0\}$ be another hypersurface such that $C$ is also self--linked
through $H$ and $G$. Recall that we defined $X:=F\cap G$. If we set
$\overline{X}:=H\cap G$, then the following proposition holds.

\proclaim{Proposition  3.1}
Let $C$ be a generically Gorenstein pure subscheme of codimension $2$ of $\p
nk$ which is self--linked through the two complete intersections
$X=F\cap G$ and
$\overline{X}=H\cap G$. If $\deg(F)\le\deg(G)$ then $X=\overline X$
\endproclaim
\demo{Proof}
At each generic point $P\in C$ we have
$\Im_{C}=(y,g)$ and $\Im_X=(f,g)$. By the factoriality of $\Ofa{\p nk,P}$
and since
$(y,g)$ is a system of parameters for the regular local ring $\Ofa{\p
nk,P}$, the
condition $\Im_X\colon\Im_C=\Im_C$ amounts to the identity $f=gz-y^2$ up to
units.

Moreover
$C$ is self--linked with respect to both $X$ and $\overline{X}$, whence in
$\Ofa{\p nk,P}$
$$
(f,g)\colon(y,g)=(y,g)=(h,g)\colon(y,g).
$$
It follows that $y^2\in(h,g)$, hence $(f,g)=(y^2,g)\subseteq(h,g)$ in
$\Ofa{\p nk,P}$. Since such an inclusion holds at each generic point $P\in C$
then
$(f,g)\subseteq(h,g)$. Changing the roles of $f$ and $h$ we obtain
$(h,g)\subseteq(f,g)$, hence equality must hold.
\qed
\enddemo

\remark{Remark 3.2}
The condition $\deg(F)\le\deg(G)$ is necessary as shows the following
easy example. Let $C$ be the origin in the affine plane with coordinates
$h,g$. Let
$f=h^2-g^2$. Then $(f,g)\ne(f,h)$.
\endremark
\medbreak

In order to clarify the role of the hypothesys that $C$ is generically
Gorenstein in
proposition 3.1, in the remaining part of this section we shall give an
example where:

\roster
\item"a)" $C$ is a pure subscheme of codimension $2$ of $\p nk$, which is not
generically Gorenstein;

\item"b)" $C$ is self--linked through $X:=F\cap G$, with $F:=\{f=0\}$ and
$G:=\{g=0\}$ of degrees $\deg(F)=:d\le m:=\deg(G)$;

\item"c)" for each choice of the second generator $\widehat{g}:=af+g$ of
$\Im_X$, there exists $h$ with $\deg(h)=\deg(f)$ such that $C$ is self--linked
through $\overline{X}:=H\cap\widehat{G}$, where
$\widehat{G}:=\{\widehat{g}=0\}$ and $H:=\{h=0\}$,

\item"d)" but $\overline{X}\ne X$.
\endroster

\example{Example 3.3}
Let $n\ge2$, $x$ and $y$ independent linear forms in $\p nk$ and $C\subseteq\p
nk$ the subscheme associated to the ideal $(x^2,xy,y^2)$. Notice that $C$ is
not generically Gorenstein.

\proclaim{Claim 3.3.1}
Let $C$ be as above. Assume that $C$ is self--linked through two hypersurfaces
$F:=\{f=0\}$ and $G:=\{g=0\}$, with $\deg(F)=:d\le
m:=\deg(G)$. Then $d=2$, $m=3$ and there exist two other linear forms
$x',y'$ on $\p nk$
such that $(x,y)=(x',y')$ and

\roster
\item"i)" either $f=x'{}^2$ and $g=y'{}^3$ mod $f$,

\item"ii)" or $f=x'y'$ and $g=x'{}^3-y'{}^3$ mod $f$.
\endroster
\endproclaim
\demo{Proof}
Since $\deg(C)=3$ then $\deg(F\cap G)=6$ hence $d=2$ and $m=3$. There exists a
Hilbert--Burch resolution  of $\Im_C$
$$
0\longrightarrow\Ofa{\p nk}(-3)^2\mapright{A}\Ofa{\p
nk}(-2)^3\mapright{B}\Im_C\longrightarrow0
$$
where
$$
A:=\pmatrix
y&0\\
-x&y\\
0&-x
\endpmatrix,\qquad B:=(x^2,xy,y^2).
$$

Then the module $\omega_C(n-4):=\Ext{1}{\p nk}\big(\Im_C,\Ofa{\p nk}(-5)\big)$
has a dual resolution
$$
0\longrightarrow\Ofa{\p nk}(-5)\mapright{{}^tB}\Ofa{\p
nk}(-3)^3\mapright{{}^tA}\Ofa{\p
nk}(-2)^2\longrightarrow\omega_C(n-4)\longrightarrow0.
$$

In particular $\omega_C(n-4)$ has two generators $e_1$, $e_2$ subject to the
relations $ye_1=xe_2=xe_1-ye_2=0$. On the other hand we have an isomorphism
$\Im_C/\Im_X\cong\omega_C(n-4)$.

Therefore $\Im_C/\Im_X$ has two generators, which lie in degree $2$,
satisfying the
above relations.

Up to a linear change $(x,y)\to(x',y')$ of generators we can assume $f=x'{}^2$
or $f=x'(x'+cy)$.

We consider then the two possible cases.

$c=0$, whence $f=x'{}^2$. Hence, modulo $f$, we can find $y'=ay+bx'$ such
that either
$g=x'y'{}^2$ or
$g=y'{}^3$. The first case is impossible since $(f,g)$ must be a regular
sequence.
In the second case let $e_1=y'{}^2,\ e_2=x'y'$ mod $\Im_X$. We obtain
therefore a
self--linkage of
$C$.

$c\ne0$ whence we may set $y'=x'+cy$ and therefore $f=x'y'$. Hence, modulo $f$,
$g=ax'{}^3+by'{}^3$ where
$ab\ne0$, thus we can assume
$g=x'{}^3-y'{}^3$. In this case $e_1=x'{}^2,\ e_2=y'{}^2$ mod $\Im_X$.\qed
\enddemo

\proclaim{Claim 3.3.2}
For each linear form $L\in k[x,y]$ there exist linear forms $\xi,\eta, M\in
k[x,y]$ such that
$$
y^3+Lx^2=\eta^3+M\xi^2
$$
where $x$ and $\xi$ (resp. $\xi$ and $\eta$) are linearly independent.
\endproclaim
\demo{Proof}
Let $L:=L_0x+L_1y$.

Assume that $L_0\ne0$. Then we may set
$$
\eta:=R'\left(x+{{L_1}\over{3L_0}}y\right),\qquad \xi=y,\qquad
M:=-{{L_1^2}\over{3L_0}}x+\left(1-{{L_1^3}\over{27L_0^2}}\right)y,
$$
where $R'$ is a fixed cube root of $L_0$.

Next consider the case $L_0=L_1=0$. Then we may take $\eta=y$, $\xi$ arbitrary
and $M=0$.

Finally let $L_0=0$ and $L_1\ne0$. Let $a$ be such that
$$
(3a^2-L_1)^2-12a^4=0.\tag3.3.2.1
$$
Notice that $a\ne0$.
It follows that $3ay^2+3a^2xy+a^3x^2-L_1xy$ is the square of a linear form
$\xi:=R''(6ay+(3a^2-L_1)x)$, where $R''$ is a fixed square root of
$(12a)^{-1}$. Setting
$\eta:=y+ax$ and
$M:=-x$, then
$y^3+Lx^2=\eta^3+M\xi^2$. Notice that $x$ and $\xi$ are independent, else $a=0$
which is not a root of equation (3.3.2.1). On the other hand also $\xi$ and
$\eta$ are independent for $L_1\ne -3a^2$, in which case substituting in
(3.3.2.1) we
would obtain $24a^4=0$, a contradiction.
\qed
\enddemo

Set now $f:=x^2$, $g:=y^3$, $\widehat{g}:=y^3+Lx^2$ and $h:=\xi^2$. Since $\xi$
and $\eta$ are independent
$(h,\widehat{g})=(\eta^3,\xi^2)$ so that $C$ is self--linked through
$H:=\{h=0\}$ and $\widehat{G}:=\{\widehat{g}=0\}$, but $f\not\in(h,g)$,
since $x$ and
$\xi$ are independent.

Thus we have checked that our example satisfies conditions a), b), c) and d)
above.
\endexample


\Refs
\refstyle{A}
\widestnumber\key{E--P--W}

\ref
\key Ap
\by R. Ap\'ery
\paper Sur certaines vari\'et\'es alg\'ebriques a $(n-2)$dimensions de
l'espace \'a
$n$ dimensions
\jour C. R. Acad. Sci. Paris
\vol 222
\yr 1946
\pages pp. 778--780
\endref

\ref
\key C-C
\by G. Casnati, F. Catanese
\paper Even sets of nodes are bundle symmetric
\jour J. Differential Geom.
\vol 47
\yr 1997
\pages pp. 237--256
\endref

\ref
\key Ca1
\by F. Catanese
\paper Babbage's conjecture, contact of surfaces, symmetric
determinantal varieties and applications
\jour Invent. Math.
\vol 63
\yr 1981
\pages pp. 433--465
\endref

\ref
\key Ca2
\by F. Catanese
\paper Homological algebra and algebraic surfaces
\inbook Algebraic Geometry, Santa Cruz 1995
\eds J. Koll\'ar, R. Lazarsfeld, D.R. Morrison
\bookinfo Proceedings of symposia in pure mathematics
\vol 62
\yr 1997
\pages pp. 3--56
\endref

\ref
\key Ei
\by D. Eisenbud
\book Commutative Algebra with a view towards Algebraic Geometry
\publ Springer
\yr 1994
\endref

\ref
\key E--P--W
\by D. Eisenbud, S. Popescu, C. Walter
\paper Lagrangian subbundles and codimension 3 subcanonical subschemes
\paperinfo Duke preprint AG/9906170
\yr 1999
\endref

\ref
\key Gae
\by F. Gaeta
\paper Sulle curve algebriche di residuale finito
\jour Ann. Mat. Pura Appl.
\vol 27
\yr 1948
\pages pp. 177--241
\endref

\ref
\key Gal
\by D. Gallarati
\paper Ricerche sul contatto di superficie algebriche lungo
curve
\jour Acad. Royale de Belgique Memoires Coll.
\vol 32
\yr 1960
\endref

\ref
\key Ha
\by R. Hartshorne
\book Algebraic geometry
\publ Springer
\yr 1977
\endref

\ref
\key P--S
\by C. Peskine, L. Szpiro
\paper Liaison des varietes algebriques
\jour Inv. Math.
\vol 26
\yr 1974
\pages 271--302
\endref

\ref
\key Rao1
\by P. Rao
\paper Liaison among curves in ${\Bbb P}^3$
\jour Inv. Math.
\vol 50
\yr 1979
\pages 205--217
\endref

\ref
\key Rao2
\by P. Rao
\paper On self--linked curves
\jour Duke Math. J.
\vol 49
\yr 1982
\pages 251--273
\endref

\ref
\key To1
\by E. Togliatti
\paper Una notevole superficie di $5^\circ$ ordine con soli punti doppi isolati
\jour Vierteljschr. Naturforsch. Ges. Z\"urich
\vol 85
\yr 1940
\pages 127--132
\endref

\ref
\key To2
\by E. Togliatti
\paper Sulle superficie col massimo numero di punti doppi
\jour Rend. Sem. Mat. Univ. Pol. Torino
\vol 9
\yr 1950
\pages 47--59
\endref

\ref
\key Wa
\by C. Walter
\paper Pfaffian subschemes
\jour J. Algebraic Geom.
\vol 5
\yr 1996
\pages pp. 671--704
\endref


\endRefs
\enddocument